\documentclass[11pt,a4paper]{amsart}
\usepackage{amsfonts}
\usepackage{amsthm}
\usepackage{amsmath}
\usepackage{amscd}
\usepackage{amssymb}
\usepackage{t1enc}
\usepackage[mathscr]{eucal}
\usepackage{indentfirst}
\usepackage{graphicx}
\usepackage{graphics}
\usepackage{pict2e}
\usepackage{epic}
\usepackage{tikz-cd}

\usepackage[style = alphabetic]{biblatex}
\renewbibmacro{in:}{}
\numberwithin{equation}{section}
\usepackage[margin=2.9cm]{geometry}
\usepackage{epstopdf}
\usepackage{comment}
\usepackage[colorinlistoftodos]{todonotes}
\usepackage[toc,page]{appendix}
\usepackage{mathtools}
\usepackage{physics}

 \definecolor{winered}{rgb}{0.5,0,0}
 \usepackage{hyperref}
 \hypersetup{citecolor = winered, linkcolor = winered, menucolor = winered, colorlinks = true}

\newtheorem{Thm}{Theorem}[section]

\theoremstyle{definition}
\newtheorem{Def}[Thm]{Definition}

\newtheorem{Ex}[Thm]{Example}

\newcommand{\cat}{\mathsf}

\RequirePackage{tikz-cd}
\RequirePackage{amssymb}
\usetikzlibrary{calc}
\usetikzlibrary{decorations.pathmorphing}

\tikzset{curve/.style={settings={#1},to path={(\tikztostart)
    .. controls ($(\tikztostart)!\pv{pos}!(\tikztotarget)!\pv{height}!270:(\tikztotarget)$)
    and ($(\tikztostart)!1-\pv{pos}!(\tikztotarget)!\pv{height}!270:(\tikztotarget)$)
    .. (\tikztotarget)\tikztonodes}},
    settings/.code={\tikzset{quiver/.cd,#1}
        \def\pv##1{\pgfkeysvalueof{/tikz/quiver/##1}}},
    quiver/.cd,pos/.initial=0.35,height/.initial=0}

\tikzset{tail reversed/.code={\pgfsetarrowsstart{tikzcd to}}}
\tikzset{2tail/.code={\pgfsetarrowsstart{Implies[reversed]}}}
\tikzset{2tail reversed/.code={\pgfsetarrowsstart{Implies}}}
\tikzset{no body/.style={/tikz/dash pattern=on 0 off 1mm}}

\author[E. Minichiello]{Emilio Minichiello}
\address{E.M., Department of Mathematics, CUNY Graduate Center}
\email{\href{mailto:eminichiello@gradcenter.cuny.edu}{eminichiello@gradcenter.cuny.edu}}

\author[M. Rivera]{Manuel Rivera}
\address{M.R., Department of Mathematics, Purdue University}
\email{\href{mailto:manuelr@purdue.edu}{manuelr@purdue.edu}}

\author[M. Zeinalian]{Mahmoud Zeinalian}
\address{M.Z., Department of Mathematics, Lehman College of CUNY}
\email{\href{mailto:mahmoud.zeinalian@lehman.cuny.edu}{mahmoud.zeinalian@lehman.cuny.edu}}

\begin{document}

\title{A detailed look at the Szczarba map}

\maketitle

\section{Introduction}
The purpose of this note is to give a detailed explanation on how to arrive to the formulae (3.11 - 3.12) in \cite{minichiello2022categorical} giving rise to an explicit natural transformation 
\[Sz \colon \mathfrak{C} \xRightarrow{ } G
\]
between two functors that we now recall. Denote by $\mathsf{sSet}$ the category of simplicial sets and by $\mathsf{Cat}_{\mathsf{sSet}}$ the category of categories enriched over the monoidal category of simplicial sets with cartesian product. The functor $\mathfrak{C} \colon \mathsf{sSet} \to \mathsf{Cat}_{\mathsf{sSet}}$ is the left adjoint of the homotopy coherent nerve functor defined by Cordier to make sense of homotopy coherent diagrams and used by Lurie to compare two models for the theory of $\infty$-categories. Conceptually, $\mathfrak{C}$ provides a combinatorial construction that allows to think of a simplicial set as a category with higher structure on the morphisms described in terms of sequences of simplices connecting two vertices. The functor $G\colon \mathsf{sSet} \to \mathsf{Cat}_{\mathsf{sSet}}$ is the left adjoint of Kan's classifying space construction usually denoted by $\overline{W} \colon \mathsf{Cat}_{\mathsf{sSet}} \to \mathsf{sSet}$, and was used by Dwyer and Kan to recast the theory of homotopy types in terms of simplicial groupoids (after formally inverting all morphisms). Conceptually, the localized version of $G$ was used in the literature to provide a combinatorial model of the path category (the many object version of the topological monoid of based loops) in terms of a simplicial set presenting a homotopy type. We recall the precise definitions of $\mathfrak{C}$ and $G$ below. 

The explicit formula we shall derive for the natural transformation $Sz$ is reminiscent of a map defined by Szczarba in terms of certain "simplicial operators" when comparing two chain models for a fibration. \cite{szczarba1961homology}. This is the reason why we call $Sz$ the \textit{Szczarba map}. 

This note may be taken as a companion to section 3 of \cite{minichiello2022categorical}, where we show that $Sz$ induces a weak equivalence after localizing both $\mathfrak{C}$ and $G$ in the context of comparing different combinatorial models for the based loop space and path category. 

\section{The functors $G$ and $\mathfrak{C}$}
Denote by $\Delta^n \in \mathsf{sSet}$ the standard $n$-simplex. Let us first define two poset-enriched categories $P_G(\Delta^n)$ and $P_{\mathfrak{C}}(\Delta^n)$. 

\begin{Def}
 The objects of $P_G(\Delta^n)$ are the elements of $[n]=\{0,1,\dots,n\}$. Given $0 \leq k \leq n$, let $\langle g_k \rangle^{n-k}$ denote the poset 
\begin{equation*}
    \langle g_k \rangle^{n-k} = \{ d_1^{n-k} g_k \leq d_1^{n-k -1} d_0 g_k \leq \dots \leq d_1 d_0^{n-k-1} g_k \leq d_0^{n-k} g_k \}.
\end{equation*}
For $n = k$, let $\langle g_k \rangle^0 = \{ g_k \}$ be the trivial poset.
The idea here is that $g_k$ is the unique non-degenerate $(n-k)$-simplex in $N\langle g_k \rangle^{n-k} \cong \Delta^{n-k}$, and the elements of the poset are the vertices of $g_k$. Recall an $\ell$-simplex $x$ in the nerve of a poset $P$ is precisely an ordered collection $(x_1 \leq x_2 \leq \dots \leq x_\ell)$ of elements in $P$.

For $0 \leq p < q \leq n$, define 
\begin{equation}
    P_G(\Delta^n)(p,q) = \langle g_q \rangle^{n-q} \times \langle g_{q-1} \rangle^{n - (q-1)} \times \dots \times \langle g_{p+1} \rangle^{n - (p + 1)},
\end{equation}
and for $p = q$, set $P_G(\Delta^n)(q,q) = \{ \text{id}_q \}$, the trivial poset. The composition rule for $P_G(\Delta^n)$ is given by the cartesian product of sequences (i.e. freely concatenating sequences). 

For example, the poset-enriched category $P_G(\Delta^2)$ looks like
\begin{equation}
\begin{tikzcd}
	&& 1 \\
	\\
	0 &&&& 2
	\arrow["{\langle g_1 \rangle^1 = \{ d_1 g_1 \leq d_0 g_1 \}}"{pos=0.3}, from=3-1, to=1-3]
	\arrow["{\langle g_2 \rangle^0 = \{ g_2 \}}", from=1-3, to=3-5]
	\arrow["{\langle g_2 \rangle^0 \times \langle g_1 \rangle^1 = \{(g_2, d_1 g_1) \leq (g_2, d_0 g_1) \}}"', shift right=2, from=3-1, to=3-5]
\end{tikzcd}
\end{equation}
Note that every arrow in the above picture is a poset.

With the cosimplicial maps given in \cite[Definition 3.2]{minichiello2022categorical}, $[n] \mapsto P_G(\Delta^n)$ defines a cosimplicial poset-enriched category.
\end{Def}

\begin{Def}
The objects of $P_{\mathfrak{C}}(\Delta^n)$ are the elements of $[n]=\{0,1,\dots,n\}$. Given $0 \leq p < q \leq n$, let $P_{\mathfrak{C}}(p,q)$ denote the poset whose elements are subsets of the form $U=\{p, i_1, \dots, i_{\ell}, q \} \subseteq \{p, p+1, \dots, q-1,q\}$ and we let $U \leq V$ if and only if $V \subseteq U$. For $p = q$, let $P_{\mathfrak{C}}(q,q) = \{ q \}$ be the trivial poset. The composition rule for $P_{\mathfrak{C}}(\Delta^n)$ is given by taking union of subsets. 

For example, the poset-enriched category $P_\mathfrak{C}(\Delta^2)$ looks like
\begin{equation*}
    \begin{tikzcd}
	&& 1 \\
	\\
	0 &&&& 2
	\arrow["{\{0, 1\}}"{pos=0.3}, from=3-1, to=1-3]
	\arrow["{\{1,2\}}", from=1-3, to=3-5]
	\arrow["{\{0,1,2\} \leq \{0,2\}}"', shift right=2, from=3-1, to=3-5]
\end{tikzcd}
\end{equation*}
Note that every arrow in the above picture is a poset.

With the cosimplicial maps given in \cite[Definition 3.7]{minichiello2022categorical}, $[n] \mapsto P_{\mathfrak{C}}(\Delta^n)$ defines a cosimplicial poset-enriched category.
\end{Def}

Any poset-enriched category gives rise to a simplicial category by applying the nerve functor at the level of posets of morphisms; we denote this functor by \[N^{\cat{Cat}} \colon \mathsf{Cat}_{\mathsf{Poset}} \to \mathsf{Cat}_{\mathsf{sSet}},\]
and define $\mathfrak{C}(\Delta^n)= N^{\cat{Cat}}( P_\mathfrak{C}(\Delta^n))$ and
$G(\Delta^n)=N^{\cat{Cat}}( P_G(\Delta^n))$. 
These constructions completely determine two functors $\mathfrak{C}$ and $G$ from $\mathsf{sSet}$ to $\mathsf{Cat}_{\mathsf{sSet}}$ through Kan extension since the category $\mathsf{Cat}_{\mathsf{sSet}}$ is cocomplete. Namely, for any arbitrary simplicial set $X$ we define
\[\mathfrak{C}(X):=\underset{{\Delta^n \to X}}{\text{colim }} \mathfrak{C}(\Delta^n)\]
and
\[G(X):=\underset{{\Delta^n \to X}}{\text{colim }}G(\Delta^n).\]

\section{An explicit formula for the Szczarba map}

In \cite[Section 2.6.1]{hinich2007homotopy}, Hinich defines a map of cosimplicial poset-enriched categories 
\[
\text{Hin} : P_\mathfrak{C}(\Delta^\bullet) \to P_G(\Delta^\bullet).
\]
It is defined as the identity on objects and on posets of morphisms 
\[\text{Hin} : P_{\mathfrak{C}}(\Delta^n)(p,q) \to P_G(\Delta^n)(p,q)\] is given on ``indecomposable" elements  by
\begin{equation*}
   \text{Hin}(\{p,q\}) = (d_1^{n-q} g_q, \, d_1^{n-q} d_0 g_{q-1}, \, \dots \, , \, d_1^{n-q} d_0^{q - (p + 1)} g_{p + 1}).
\end{equation*}
Then on any element $U=\{p, i_1, \dots, i_\ell, q \} \in P_{\mathfrak{C}}(\Delta^n)(p,q)$, with $p < i_1 < i_2 < \dots < i_\ell <q$, one may define
\begin{equation*}
    \text{Hin}(\{p, i_1, \dots, i_\ell, q \}) = (\text{Hin}(\{i_\ell, q \}), \, \dots \, , \, \text{Hin}(\{ p, i_1 \}) ).
\end{equation*}
 
\begin{Ex}
For $\{0, 2, 4 \}$ in $P_{\mathfrak{C}}(\Delta^5)(0,5)$, we have
\begin{equation*}
    \begin{aligned}
   \text{Hin}(\{0, 2, 4\}) & = (\text{Hin}(\{2,  4\}), \text{Hin}(\{0,2 \}))  \\
   & = (d_1 g_4, d_1 d_0 g_3 , d_1^3 g_2, d_1^3 d_0 g_1). 
    \end{aligned}
\end{equation*}
\end{Ex}

The goal is to obtain an explicit formula for the morphism of simplicial categories
\[N^{\cat{Cat}} \text{Hin} : \mathfrak{C}(\Delta^n) \to G(\Delta^n).\] We shall denote this functor of simplicial categories by $Sz_{\Delta^n}$.

It is enough to consider the non-degenerate $\ell$-simplices of $\mathfrak{C}(\Delta^n)(p,q)$. For a fixed $n \geq 1$, let $0 \leq p < q \leq n$ and $0 \leq \ell \leq q-p-1$. Consider the set $\text{nd}(\mathfrak{C}(\Delta^n)(p,q)_\ell)$ of non-degenerate $\ell$-simplices of $\mathfrak{C}(\Delta^n)(p,q)$. There is an obvious bijection
$\text{nd}(\mathfrak{C}(\Delta^n)(p,q)_\ell) \cong S^\ell_{p,q}$, where
\begin{equation*}
    S^\ell_{p,q} = \{ i = (i_1, \dots, i_\ell) \in \{ p + 1, \dots, q -1 \}^{\times \ell} \; : \: i_r \neq i_s, \text{ for } r \neq s \}.
\end{equation*}
For $\ell = 0$, we set $S^0_{p,q} = \{ \emptyset \}$, and we call $\emptyset$ the empty sequence. From now on, we identify sequences $(i_1, \dots, i_\ell)$ as above with non-degenerate $\ell$-simplices in $\mathfrak{C}(\Delta^n)(p,q)$.

Given an $\ell$-simplex
$$i = (\{p, i_1, \dots, i_\ell, q \} \leq \{p, i_1, \dots, i_{\ell - 1}, q \} \leq \dots \leq \{p, i_1, q \} \leq \{p, q \}),$$
we shall describe a formula for the resulting $\ell$-simplex \[Sz_{\Delta^n}(i) = N^{\cat{Cat}} \text{Hin}(i) \in G(\Delta^n)(p,q)_\ell = N(\langle g_q \rangle^{n-q})_{\ell} \times N(\langle g_{q-1} \rangle^{n - (q-1)})_{\ell} \times \dots \times N(\langle g_{p+1} \rangle^{n - (p + 1)})_{\ell}.\]

Therefore $Sz_{\Delta^n}(i)$ will have components corresponding to each $g_k$ for $p + 1 \leq k \leq q$. For each fixed $k$, we shall give a formula for the $k$-th component $Sz_{\Delta^n}(i)_k = \mathcal{E}_{i,k} g_k$, in terms of a simplicial operator $\mathcal{E}_{i,k}$ that will be defined inductively.  Note $\mathcal{E}_{i,k} g_k$ must be an $\ell$-simplex in $N (\langle g_k \rangle^{n-k})\cong \Delta^{n-k}$, so it is enough to specify all of its $\ell +1$ vertices. Lets use the notation $[x_\ell, x_{\ell - 1}, \dots, x_0]$, where $0\leq x_\ell \leq \dots \leq x_0 \leq n-k$, to denote the $\ell$-simplex corresponding to $\mathcal{E}_{i,k} g_k$ in $\Delta^{n-k}$. 

\begin{Ex}
For the $2$-simplex corresponding to
\begin{equation*}
   \{0,3 \} \geq \{0, 2, 3 \} \geq \{0, 1, 2, 3\} 
\end{equation*}
we obtain the $2$-simplex
\begin{equation*}
    (g_3, d_0 g_2, d_0^2 g_1) \geq (g_3, d_1 g_2, d_1 d_0 g_1) \geq (g_3, d_1 g_2, d_1^2 g_1).
\end{equation*}
Writing this component-wise as
\begin{equation*}
\begin{aligned}
    & g_3 \leq g_3 \leq g_3 \\
    & d_1 g_2 \leq d_1 g_2 \leq d_0 g_2 \\
    & d_1^2 g_1 \leq d_1 d_0 g_1 \leq d_0^2 g_1
\end{aligned}
\end{equation*}
we can see that for each $k$ with $1 \leq k \leq 3$, we obtain a $2$-simplex in $N \langle g_k \rangle^{3 - k} \cong \Delta^{3-k}$. These are given by
\begin{equation*}
    \begin{aligned}
    &k = 3, \qquad  [0 \, 0 \, 0] \\
    &k = 2, \qquad  [0 \, 0 \, 1] \\
    &k = 1, \qquad [0 \, 1 \, 2].
    \end{aligned}
\end{equation*}
Thus for $k = 2$, we have $x_0 = 1$, $x_1 = 0$ and $x_2 = 0$.
\end{Ex}

For a fixed $k$, it is easy to obtain $x_0$, as it will appear in $\text{Hin}(\{p, q \})$ as $d_1^{n-q} d_0^{q-k} g_k$, which corresponds to the vertex $q-k$ in the poset $\langle g_k \rangle^{n-k}$. Thus $x_0 = q - k$ for every $k$. 

Now suppose we start with $i = \emptyset$, corresponding to $\{p, q\}$ and $x_0$ and add $i_1$ to it. We thus wish to compute $x_1$. But we will need to know the data of $i_1$. In fact all the higher $x_\ell$ will depend on the sequence $(i_1, \dots, i_{\ell - 1})$. We define a function
\[\alpha_k : S^\ell_{p,q} \to \{0, \dots, n - k\}\]
so that
\begin{equation*}
    x_\ell = \alpha_k(i_1, \dots, i_{\ell - 1}).
\end{equation*}
 As mentioned, we have $\alpha_k(\emptyset) = q - k$.

When we add $i_1$, we obtain a $1$-simplex $\text{Hin}(\{p, q\}) \geq \text{Hin}(\{p,i_1, q \})$, which gives

\begin{equation*}
\begin{aligned}
    & (d_1^{n-q} g_q, \dots, d_1^{n-q} d_0^{q - (i_1 + 1)} g_{i_1 + 1}, d_1^{n-q} d_0^{q - i_1} g_{i_1}, d_1^{n-q} d_0^{q - (i_1 - 1)} g_{i_1 - 1}, \dots, d_1^{n-q} d_0^{q - (p + 1)} g_{p+1}) \\
   &  \geq (d_1^{n-q} g_q, \dots, d_1^{n-q} d_0^{q - (i_1 + 1)} g_{i_1 + 1}, \, d_1^{n - i_1} g_{i_1}, \, d_1^{n - i_1} d_0 g_{i_1 - 1}, \, \dots \,, d_1^{n - i_1} d_0^{i_1 - (p + 1)} g_{p+1})
\end{aligned}
\end{equation*}
Note that for $k > i_1$, the coefficients for $g_k$ do not change. For $k \leq i_1$, we see that the vertex $x_1$ corresponds to $d_1^{n - i_1} d_0^{i_1 - k} g_k$, which is the $(i_1 - k)$-th element of the linear ordering $\langle g_k \rangle^{n-k}$. Thus we have
\begin{equation*}
    x_1 = \alpha_k(i_1) = \begin{cases}
        \alpha_k(\emptyset) & k > i_1 \\
        i_1 - k & k \leq i_1.
    \end{cases}
\end{equation*}

We now compute $\alpha_k(i_1, \dots, i_{\ell})$ given we have computed \[\alpha_k(\emptyset), \alpha_k(i_1), \dots, \alpha_k(i_1, \dots, i_{\ell - 1}).\]

This corresponds to the $1$-simplex $\text{Hin}(\{p, i_1, \dots, i_{\ell - 1}, q\}) \geq \text{Hin}(\{p, i_1, \dots, i_{\ell}, q \})$. In order to see which of the $g_k$ will be affected, we must know where $i_\ell$ appears amongst the $i_1, \dots, i_{\ell - 1}$ if they are put in order.

For this we define $\omega_{(i_1, \dots, i_{\ell - 1})}(i_\ell)$ to be the largest integer in $\{p, i_1, \dots, i_{\ell - 1}, q \}$ such that $\omega_{(i_1, \dots, i_{\ell - 1})}(i_{\ell}) < i_\ell$. This number will tell us the lower bound for the range of $k$ for which the coefficients of $g_k$ will not change. Within the range $\omega_{(i_1, \dots, i_{\ell - 1})}(i_{\ell}) < k \leq i_{\ell}$, we know that the coefficient of $g_k$ will change to $d_1^{n- i_{\ell}} d_0^{i_{\ell} - k} g_k$, by inspecting what happens within Hinich's formula. This is the $(i_\ell - k)$-th element of the linear ordering of $\langle g_k \rangle^{n- k}$. Outside of this range, we know that the coefficients for $g_k$ will not change. Thus we define
\begin{equation}
   \alpha_k(i_1, \dots, i_{\ell}) = 
   \begin{cases}
  i_\ell - k, & \omega_{(i_1, \dots, i_{\ell - 1})}(i_\ell) < k \leq i_\ell \\
  \alpha_k(i_1, \dots, i_{\ell - 1}), & k \leq \omega_{(i_1, \dots, i_{\ell - 1})}(i_\ell) \text{ or } i_\ell < k.
    \end{cases}
\end{equation}
For any sequence $i = (i_1, \dots, i_\ell)$ we have constructed $[x_\ell \; x_{\ell - 1} \dots x_0]$ corresponding to the $\ell$-simplex $\mathcal{E}_{i,k} g_k$ in $N (\langle g_k \rangle^{n-k}) \cong \Delta^{n-k}$. From the sequence $[x_\ell \; x_{\ell - 1} \dots x_0]$ we wish to obtain a simplicial operator $\mathcal{E}_{i,k}$, defined inductively, such that $x_r = d_1^{\ell - r} d_0^r (\mathcal{E}_{i,k} g_k)$. Namely, $x_r$ is the $r$th vertex of the $\ell$-simplex $\mathcal{E}_{i,k} g_k$.

So assuming $\mathcal{E}_{i', k} g_k = [x_{\ell - 1} \dots x_0]$, where $i' = (i_1, \dots, i_{\ell - 1})$, we wish to know what simplicial operator we need in order to produce $\mathcal{E}_{i,k}g_k = [x_\ell \; x_{\ell - 1} \dots x_0]$ thinking of $g_k$ as $[0\;1\; \dots \;(n-k)]$. If $x_{\ell} = x_{\ell - 1}$, then this is easy, namely we need only set $\mathcal{E}_{i,k} = s_0 \mathcal{E}_{i', k}$, since 
\begin{equation*}
    [x_{\ell - 1} \; x_{\ell - 1} \dots x_0] = s_0 [x_{\ell - 1} \dots x_0]. 
\end{equation*}
If $x_\ell < x_{\ell - 1}$, we can obtain $\mathcal{E}_{i,k} g_k$ from $\mathcal{E}_{i', k} g_k$ by noticing 
\begin{equation*}
     s_0^{x_{\ell} + 1} d_0^{x_\ell} [0\; 1 \dots (n - k)]= 
     [x_{\ell} \dots x_{\ell} \dots x_{\ell} \; (x_{\ell} + 1) \dots x_{\ell-1} \dots x_{\ell-2} \dots x_0 \dots (n - k)]
\end{equation*}
Then applying $\mathcal{E}'_{i',k}$, where the superscript $'$ means to add one to every index in the simplicial operator, will "skip over" the first $x_\ell$ and will then "cross out" everything else but $\mathcal{E}_{i',k} g_k = [x_{\ell - 1} \; x_{\ell - 2} \; \dots x_0]$, but now this will be placed after the additional $x_{\ell}$, leaving
\begin{equation*}
    \mathcal{E}_{i,k} g_k = [x_{\ell} \; x_{\ell - 1} \dots x_0].
\end{equation*}
It follows that the desired operators $\mathcal{E}_{i,k}$ may be defined by induction on $\ell$, the length of $i$, as follows. For the empty sequence, define
\begin{equation*}
        \mathcal{E}_{\emptyset, k} = d_1^{n-q} \, d_0^{q-k}.
\end{equation*}
Suppose we have defined $\mathcal{E}_{i,k}$ for any $i$ of length less than $\ell$.  If $i$ has length $\ell$, define
\begin{equation}
    \mathcal{E}_{i,k} = 
    \begin{cases}
        s_0 \, \mathcal{E}_{i',k} & \text{if } \alpha_k(i') = \alpha_k(i) \\
        \mathcal{E}'_{i', k} \, s_0^{\alpha_{k}(i) + 1} \, d_0^{\alpha_k(i)} & \text{if } \alpha_k(i) < \alpha_k(i').
    \end{cases}
\end{equation}
In summary, we have shown that the map \[Sz_{\Delta^n}= N^{\cat{Cat}}(\text{Hin}) \colon  \mathfrak{C}(\Delta^n)(p,q) \to G(\Delta^n)(p,q)\] is given on any 
\[i = (\{p, i_1, \dots, i_\ell, q \} \leq \{p, i_1, \dots, i_{\ell - 1}, q \} \leq \dots \leq \{p, i_1, q \} \leq \{p, q \}) \in \mathfrak{C}(\Delta^n)(p,q)_{\ell}\]
by the formula
\begin{equation*} \label{Sz formula}
Sz_{\Delta^n}(i) =  \left( \mathcal{E}_{i,q} \, g_q, \, \mathcal{E}_{i,q-1} \, g_{q-1}, \, \dots \, , \, \mathcal{E}_{i, p+2} \, g_{p + 2}, \, \mathcal{E}_{i, p+1} \, g_{p+1} \right) \in G(\Delta^n)(p,q)_{\ell}.
\end{equation*}

We finish by repeating Example 3.15 from \cite{minichiello2022categorical}.

\begin{Ex}
Consider the $2$-simplex in $\mathfrak{C}(\Delta^3)(0,3)_2$ given by 
$$
\{0,3\} \geq \{0,2,3\} \geq \{0,1,2,3\}.
$$
This simplex corresponds to the sequence $i = (2,1)$. So with $n = 3, \, p = 0, \, q = 3$, we compute
\begin{equation*}
    \alpha_3(\emptyset) = 0, \; \alpha_2(\emptyset) = 1, \; \alpha_1(\emptyset) = 2
\end{equation*}
\begin{equation*}
   \omega_\emptyset(2) = 0, \; \alpha_3(2) = 0, \; \alpha_2(2) = 0, \; \alpha_1(2) = 1 
\end{equation*}
\begin{equation}
\omega_{(2)}(1) = 0, \; \alpha_3(2,1) = 0, \; \alpha_2(2,1) = 0, \; \alpha_1(2,1) = 0.
\end{equation}
With this we can then compute
\begin{equation}
    \begin{aligned}
    Sz_{\Delta^3}(i) & = (\mathcal{E}_{(2,1), 3} \, g_3, \, \mathcal{E}_{(2,1),2} \, g_2, \, \mathcal{E}_{(2,1),1} \, g_1) \\
    & = (s_0 \mathcal{E}_{(2),3} \, g_3, \, s_0 \mathcal{E}_{(2),2} \, g_2, \, \mathcal{E}'_{(2), 1} s_0 g_1) \\
    & = (s_0^2 \mathcal{E}_{\emptyset, 3} \, g_3, \, s_0 \mathcal{E}'_{\emptyset, 2} s_0 g_2, \, \mathcal{E}''_{\emptyset, 1} s_1^2 d_1 s_0 g_1)  \\
    & = (s_0^2 g_3, \, s_0 d_1 s_0 g_2, \, d_2^2 s_1^2 g_1) \\
    & = (s_0^2 g_3, \, s_0 g_2, \, g_1).
    \end{aligned}
\end{equation} 
A similar computation gives 
\begin{equation*}
    Sz_{\Delta^3}(1,2) = (s_0^2 g_3, \, s_1 g_2, \, s_0 d_1 g_1).
\end{equation*}
\end{Ex}

We include a diagram illustrating the map $Sz_{\Delta^3} : \mathfrak{C}(\Delta^3)(0,3) \to G(\Delta^3)(0,3)$.

\vspace{.5cm}

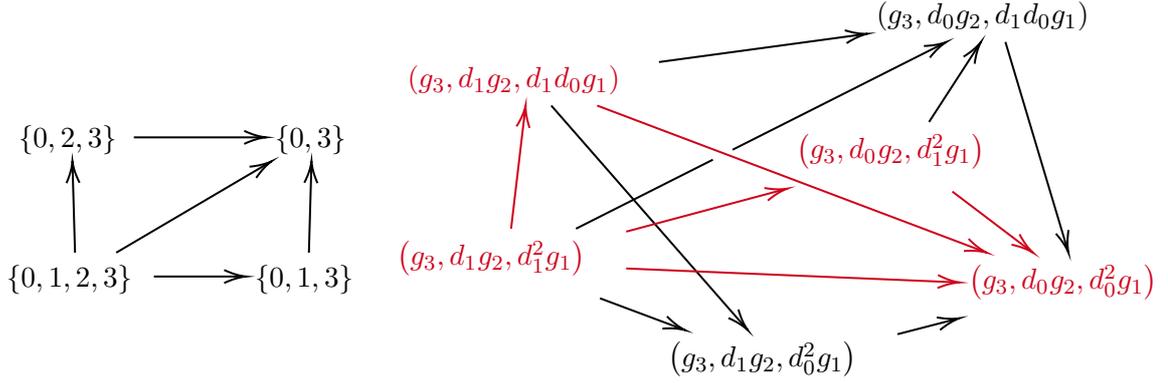
\tikzset{every picture/.style={line width=0.75pt}} 

\begin{tikzpicture}[x=0.75pt,y=0.75pt,yscale=-1,xscale=1]

\draw (155,132.4) node [anchor=north west][inner sep=0.75pt]    {$\{0,3\}$};
\draw (145,202.4) node [anchor=north west][inner sep=0.75pt]    {$\{0,1,3\}$};
\draw (21,202.4) node [anchor=north west][inner sep=0.75pt]    {$\{0,1,2,3\}$};
\draw (27,132.4) node [anchor=north west][inner sep=0.75pt]    {$\{0,2,3\}$};
\draw (215.67,190.35) node [anchor=north west][inner sep=0.75pt]  [color={rgb, 255:red, 208; green, 2; blue, 27 }  ,opacity=1 ]  {$\left( g_{3} ,d_{1} g_{2} ,d_{1}^{2} g_{1}\right)$};
\draw (221,102.4) node [anchor=north west][inner sep=0.75pt]  [color={rgb, 255:red, 208; green, 2; blue, 27 }  ,opacity=1 ]  {$( g_{3} ,d_{1} g_{2} ,d_{1} d_{0} g_{1})$};
\draw (351,241.4) node [anchor=north west][inner sep=0.75pt]    {$\left( g_{3} ,d_{1} g_{2} ,d_{0}^{2} g_{1}\right)$};
\draw (455,70.4) node [anchor=north west][inner sep=0.75pt]    {$( g_{3} ,d_{0} g_{2} ,d_{1} d_{0} g_{1})$};
\draw (415,136.66) node [anchor=north west][inner sep=0.75pt]  [color={rgb, 255:red, 208; green, 2; blue, 27 }  ,opacity=1 ]  {$\left( g_{3} ,d_{0} g_{2} ,d_{1}^{2} g_{1}\right)$};
\draw (501,202.4) node [anchor=north west][inner sep=0.75pt]  [color={rgb, 255:red, 208; green, 2; blue, 27 }  ,opacity=1 ]  {$\left( g_{3} ,d_{0} g_{2} ,d_{0}^{2} g_{1}\right)$};
\draw    (96,210) -- (140,210) ;
\draw [shift={(142,210)}, rotate = 180] [color={rgb, 255:red, 0; green, 0; blue, 0 }  ][line width=0.75]    (10.93,-3.29) .. controls (6.95,-1.4) and (3.31,-0.3) .. (0,0) .. controls (3.31,0.3) and (6.95,1.4) .. (10.93,3.29)   ;
\draw    (173.34,198) -- (174.6,154) ;
\draw [shift={(174.66,152)}, rotate = 91.64] [color={rgb, 255:red, 0; green, 0; blue, 0 }  ][line width=0.75]    (10.93,-3.29) .. controls (6.95,-1.4) and (3.31,-0.3) .. (0,0) .. controls (3.31,0.3) and (6.95,1.4) .. (10.93,3.29)   ;
\draw    (56.66,198) -- (55.4,154) ;
\draw [shift={(55.34,152)}, rotate = 88.36] [color={rgb, 255:red, 0; green, 0; blue, 0 }  ][line width=0.75]    (10.93,-3.29) .. controls (6.95,-1.4) and (3.31,-0.3) .. (0,0) .. controls (3.31,0.3) and (6.95,1.4) .. (10.93,3.29)   ;
\draw    (86,140) -- (150,140) ;
\draw [shift={(152,140)}, rotate = 180] [color={rgb, 255:red, 0; green, 0; blue, 0 }  ][line width=0.75]    (10.93,-3.29) .. controls (6.95,-1.4) and (3.31,-0.3) .. (0,0) .. controls (3.31,0.3) and (6.95,1.4) .. (10.93,3.29)   ;
\draw    (77.23,198) -- (153.05,153.02) ;
\draw [shift={(154.77,152)}, rotate = 149.32] [color={rgb, 255:red, 0; green, 0; blue, 0 }  ][line width=0.75]    (10.93,-3.29) .. controls (6.95,-1.4) and (3.31,-0.3) .. (0,0) .. controls (3.31,0.3) and (6.95,1.4) .. (10.93,3.29)   ;
\draw [color={rgb, 255:red, 208; green, 2; blue, 27 }  ,draw opacity=1 ]   (274.22,185.95) -- (281.24,125.99) ;
\draw [shift={(281.48,124)}, rotate = 96.68] [color={rgb, 255:red, 208; green, 2; blue, 27 }  ,draw opacity=1 ][line width=0.75]    (10.93,-3.29) .. controls (6.95,-1.4) and (3.31,-0.3) .. (0,0) .. controls (3.31,0.3) and (6.95,1.4) .. (10.93,3.29)   ;
\draw    (318.56,220.95) -- (359.24,236.29) ;
\draw [shift={(361.11,237)}, rotate = 200.67] [color={rgb, 255:red, 0; green, 0; blue, 0 }  ][line width=0.75]    (10.93,-3.29) .. controls (6.95,-1.4) and (3.31,-0.3) .. (0,0) .. controls (3.31,0.3) and (6.95,1.4) .. (10.93,3.29)   ;
\draw    (294.28,124) -- (391.01,235.49) ;
\draw [shift={(392.32,237)}, rotate = 229.06] [color={rgb, 255:red, 0; green, 0; blue, 0 }  ][line width=0.75]    (10.93,-3.29) .. controls (6.95,-1.4) and (3.31,-0.3) .. (0,0) .. controls (3.31,0.3) and (6.95,1.4) .. (10.93,3.29)   ;
\draw    (348,102.11) -- (450.02,88.16) ;
\draw [shift={(452,87.89)}, rotate = 172.21] [color={rgb, 255:red, 0; green, 0; blue, 0 }  ][line width=0.75]    (10.93,-3.29) .. controls (6.95,-1.4) and (3.31,-0.3) .. (0,0) .. controls (3.31,0.3) and (6.95,1.4) .. (10.93,3.29)   ;
\draw [color={rgb, 255:red, 208; green, 2; blue, 27 }  ,draw opacity=1 ]   (331.67,187.42) -- (410.07,166.31) ;
\draw [shift={(412,165.79)}, rotate = 164.93] [color={rgb, 255:red, 208; green, 2; blue, 27 }  ,draw opacity=1 ][line width=0.75]    (10.93,-3.29) .. controls (6.95,-1.4) and (3.31,-0.3) .. (0,0) .. controls (3.31,0.3) and (6.95,1.4) .. (10.93,3.29)   ;
\draw    (482.75,132.26) -- (507.56,93.68) ;
\draw [shift={(508.64,92)}, rotate = 122.74] [color={rgb, 255:red, 0; green, 0; blue, 0 }  ][line width=0.75]    (10.93,-3.29) .. controls (6.95,-1.4) and (3.31,-0.3) .. (0,0) .. controls (3.31,0.3) and (6.95,1.4) .. (10.93,3.29)   ;
\draw    (520.86,92) -- (551.74,196.08) ;
\draw [shift={(552.31,198)}, rotate = 253.47] [color={rgb, 255:red, 0; green, 0; blue, 0 }  ][line width=0.75]    (10.93,-3.29) .. controls (6.95,-1.4) and (3.31,-0.3) .. (0,0) .. controls (3.31,0.3) and (6.95,1.4) .. (10.93,3.29)   ;
\draw [color={rgb, 255:red, 208; green, 2; blue, 27 }  ,draw opacity=1 ]   (494.39,167.26) -- (533.02,196.79) ;
\draw [shift={(534.61,198)}, rotate = 217.39] [color={rgb, 255:red, 208; green, 2; blue, 27 }  ,draw opacity=1 ][line width=0.75]    (10.93,-3.29) .. controls (6.95,-1.4) and (3.31,-0.3) .. (0,0) .. controls (3.31,0.3) and (6.95,1.4) .. (10.93,3.29)   ;
\draw    (467,239.03) -- (496.06,231.47) ;
\draw [shift={(498,230.97)}, rotate = 165.43] [color={rgb, 255:red, 0; green, 0; blue, 0 }  ][line width=0.75]    (10.93,-3.29) .. controls (6.95,-1.4) and (3.31,-0.3) .. (0,0) .. controls (3.31,0.3) and (6.95,1.4) .. (10.93,3.29)   ;
\draw [color={rgb, 255:red, 208; green, 2; blue, 27 }  ,draw opacity=1 ]   (331.67,205.96) -- (496,212.9) ;
\draw [shift={(498,212.99)}, rotate = 182.42] [color={rgb, 255:red, 208; green, 2; blue, 27 }  ,draw opacity=1 ][line width=0.75]    (10.93,-3.29) .. controls (6.95,-1.4) and (3.31,-0.3) .. (0,0) .. controls (3.31,0.3) and (6.95,1.4) .. (10.93,3.29)   ;
\draw [color={rgb, 255:red, 208; green, 2; blue, 27 }  ,draw opacity=1 ]   (317.15,124) -- (509.66,197.29) ;
\draw [shift={(511.53,198)}, rotate = 200.84] [color={rgb, 255:red, 208; green, 2; blue, 27 }  ,draw opacity=1 ][line width=0.75]    (10.93,-3.29) .. controls (6.95,-1.4) and (3.31,-0.3) .. (0,0) .. controls (3.31,0.3) and (6.95,1.4) .. (10.93,3.29)   ;
\draw    (306.6,185.95) -- (374.91,151.22)(384.71,146.24) -- (489.64,92.91) ;
\draw [shift={(491.42,92)}, rotate = 153.06] [color={rgb, 255:red, 0; green, 0; blue, 0 }  ][line width=0.75]    (10.93,-3.29) .. controls (6.95,-1.4) and (3.31,-0.3) .. (0,0) .. controls (3.31,0.3) and (6.95,1.4) .. (10.93,3.29)   ;

\end{tikzpicture}

\vspace{.5cm}

The diagram on the left is an illustration of the nondegenerate simplices in $\mathfrak{C}(\Delta^3)(0,3) \cong \Delta^1 \times \Delta^1$ and similarly on the right for $G(\Delta^3)(0,3) \cong \Delta^0 \times \Delta^1 \times \Delta^2$. The red subdiagram on the right shows the image of the Szczarba map.

\printbibliography

\end{document}